\documentclass[12pt]{amsart}
\usepackage{amsmath,amsthm,amssymb}
\usepackage{mathrsfs}
\newtheorem{theorem}{Theorem}[section]
\newtheorem{lemma}[theorem]{Lemma}
\newtheorem{proposition}[theorem]{Proposition}

\newtheorem{definition}[theorem]{Definition}
\numberwithin{equation}{section}

\begin{document}

\title{An extension of harmonic functions along fixed direction}

\author{Sevdiyor Imomkulov,
        Yuldash Saidov}


\abstract Let a function $u(x,y)$ be harmonic in the domain
$$
D\times V_r=D\times \{y\in \mathbb{R}^m: |y|<r\}\subset
\mathbb{R}^n\times \mathbb{R}^m
$$
and for each fixed point $x^0$ from some a set $E\subset D$, 
the function $u(x^0,y)$, as a function of variable $y$, 
can be extended to a harmonic
function on the whole $\mathbb{R}^m$. Then $u(x,y)$ harmonically extends
to the domain $D\times \mathbb{R}^m$ as a function of variables $x$ and
$y$.
\endabstract

\keywords{Harmonic function; harmonic extension; Hartogs theorem;
Hartog's series; polar sets; holomorphic function; function class
$Lh_{0}(D)$.}

\subjclass{ 32A65, 46J10, 46J15.}

\maketitle

\section
{ Intoduction}

The well known Hartogs lemma concerns the extension of holomorphic functions
along fixed direction (see \cite{shabat}) and states that  if a
function $f(z,w)$ is holomorphic in the domain
$$
D\times \{w\in \mathbb{C}: |w|<r\}\subset \mathbb{C}^n\times
\mathbb{C}
$$
and for each fixed $z^0\in D$ the function $f(z^0,w)$ of variable
$w$  holomorphically extends to the disk $ \{w\in \mathbb{C}:
|w|<R\}, R>r$, then $f(z,w)$ can be  holomorphically extended to the
domain $ D\times \{w\in \mathbb{C}: |w|<R\}$ as a function of variables
$z$ and $w$.

Note that the Hartogs lemma remains to hold also for pluriharmonic and harmonic
functions (see \cite{sadullaev}) and can be formulated as follows. 
Let a function $u(x,y)$ be harmonic in the domain
$$
D\times V_r=D\times \{y\in \mathbb{R}^m: |y|<r\}\subset
\mathbb{R}^n\times \mathbb{R}^m,
$$
and for each fixed $x^0\in D$ the function $u(x^0,y)$ of variable
$y$  harmonically extends to the ball $V_R=\{y\in \mathbb{R}^m:
|y|<R\}, R>r.$ Then $u(x,y)$ harmonically extends to $D\times V_R$
as a function of variables $x$ and $y$.

One of the main methods used in the investigation of harmonic extensions
is to convert this problem to holomorfic extensions.
For this purpose, the following lemma is proved in \cite{sadullaev}. 

\begin{lemma}                                  \label{L1}
Consider the space $\mathbb{R}^n(x)$ embedded into
$\mathbb{C}^n(z)=\mathbb{R}^n(x)+i\mathbb{R}^n(y)$, where
$z=(z_1,z_2,...,z_n),~z_j=x_j+iy_j,j=1,2,...,n$. Then for each
domain $D\subset \mathbb{R}^n(x)$ there exists some holomorphy domain
$\widehat{D}\subset \mathbb{C}^n(z)$, such that $D\subset
\widehat{D}$ and each harmonic function in $D$ can be holomorphically
extended to the domain $\widehat{D}$, i.e., there exists a holomorphic
function $f_u(z)\in \widehat{D}$ such that $f_u|_D=u.$
\end{lemma}

The main result of our paper reads as follows. The
concept of $N$-sets of $  Lh_0(D)$ which appears there
is defined in the next section.

\begin{theorem}                             \label{T1}
Let a function $u(x,y)$ be harmonic in the domain
$$
D\times V_r=D\times \{y\in \mathbb{R}^m: |y|<r\}\subset
\mathbb{R}^n\times \mathbb{R}^m
$$
and for each fixed point $x^0$ from some set $E\subset D$, which is
not embedded into a countable union of $N$-sets of $Lh_0(D)$,
the function $u(x^0,y)$ of variable $y$ can be extended to a harmonic
function on the whole $\mathbb{R}^m$. Then $u(x,y)$ harmonically extends
to the domain $D\times \mathbb{R}^m$ as a function of variables $x$ and $y$.
\end{theorem}

\section {\bf The class of functions $Lh_0(D)$}

Let $D\subset \mathbb{R}^n$ and $h(D)$ be the space of harmonic
functions in $D$. By $Lh_\varepsilon(D)$ we denote the minimal class
of functions which contains functions of the form $\alpha\ln
|u(x)|,~ u(x)\in h(D),~ 0<\alpha<\varepsilon$, and it is closed
with respect to the  operation of ``upper regularization'', i.e., 
for any set of
functions $u_\lambda(x)\in Lh_\varepsilon(D),~\lambda\in \Lambda$,
where $\lambda$ is an index set, the function
$$
\overline {{\mathop {\lim} \limits_{y \to x}}}  (\sup \{u_{\lambda}
(y):\ \lambda \in \Lambda \})
$$
also belongs to the class $Lh_\varepsilon(D)$(\cite{zahariuta}).

\noindent The union of function classes
$Lh_0(D)=\bigcup\limits_{\varepsilon>0}Lh_\varepsilon(D)$ we call
{\it the set of $Lh_0-$ functions}.

Now, we will define $N$-set of the class $Lh_0(D)$. Let
$\vartheta_k(x)\in Lh_0(D)$ be a monotonically increasing sequence of
locally uniformly from above bounded functions. We denote

$$
\vartheta(x)=\overline {{\mathop {\lim}\limits_{y\rightarrow
x}}}\lim\limits_{k\rightarrow\infty}\vartheta_k(y)=
\ x\in D.
$$
Then everywhere in $D$ holds the inequality
$$
\lim\limits_{k\rightarrow\infty}\vartheta_k(y)\leq\vartheta(x).
$$

\begin{definition}                    \label{D1}
The subsets of the set
$$
\{x\in
D:\lim\limits_{k\rightarrow\infty}\vartheta_k(y)<\vartheta(x)\}
$$
are called $N$-sets of the class $Lh_0(D)$. 
\end{definition}

Note that if $\vartheta_k(x)\in Lh_0(D)$ is a sequence of locally
uniformly from above bounded functions and
$$
{\overline{\mathop {\lim}\limits_{y\rightarrow x}}}\,\,\,
{\overline{\mathop{\lim}\limits_{k\rightarrow\infty}}}
\vartheta_k(y)=\vartheta(x),~~x\in D,
$$
then the set
$$
E=\{x\in D:\overline{\mathop
{\lim}\limits_{k\rightarrow\infty}}\vartheta_k(y)<\vartheta(x)\}
$$
consists of countable union of $N$-sets of the class
$Lh_0(D)$. Indeed, consider the sequence of functions
$$
w_{l,j}(x)=\max\limits_{l\leq k\leq j}\vartheta_k(x).
$$
Clearly
$$
\overline{\mathop
{\lim}\limits_{k\rightarrow\infty}}\vartheta_k(y)=\lim
\limits_{l\rightarrow \infty}\lim \limits_{j\rightarrow\infty}
w_{l,j}(x).
$$
Since the sequence is increasing in $j$, we have
$$
\lim \limits_{j\rightarrow \infty} w_{l,j}(x)\leq \overline{\mathop
{\lim} \limits_{l\rightarrow \infty}}\lim
\limits_{j\rightarrow\infty} w_{l,j}(y),~~x\in D,
$$
and the sets
$$
E_l=\{x\in D:\lim \limits_{j\rightarrow \infty}
w_{l,j}(x)<\overline{\mathop {\lim}\limits_{l\rightarrow
\infty}}\lim \limits_{j\rightarrow\infty} w_{l,j}(y)\},~l=1,2,\dots,
$$
are $N$-set of class the $Lh_0(D)$. On the other hand, the sequences
$\lim \limits_{j\rightarrow \infty} w_{l,j}(x)$ and
$\overline{\mathop {\lim}\limits_{l\rightarrow \infty}}\lim
\limits_{j\rightarrow\infty} w_{l,j}(y),~j=1,2,\dots,$ are
monotonically decreasing, and
$$
\overline{\mathop
{\lim}\limits_{k\rightarrow\infty}}\vartheta_k(x)=\lim
\limits_{l\rightarrow \infty}\lim \limits_{j\rightarrow\infty}
w_{l,j}(x)=\vartheta(x)=\lim \limits_{l\rightarrow
\infty}\overline{\mathop {\lim}\limits_{y\rightarrow x}}\lim
\limits_{j\rightarrow\infty} w_{l,j}(y),~x\in D\setminus
\bigcup\limits_{l=1}^\infty E_l.
$$
It follows that $E\subset \bigcup\limits_{l=1}^\infty E_l$, i.e.
$E=\bigcup\limits_{l=1}^\infty (E_l\cap E)$.

\begin{definition}                             \label{D2}
The set $E\subset D$ is called $Lh_0$-{\it polar
with respect to $D$}, if there exists a function $\vartheta(x)\in
Lh_0(D)$ such that $\vartheta(x){\not \equiv}-\infty$ and
$\vartheta(x)|_E=-\infty.$ 
\end{definition}

Note that if $u(x)\in h(D),~u(x){\not \equiv}0$  and $E\subset
\{u(x)= 0\}$, then $E$ is $Lh_0$-polar with respect to $D$.

\begin{proposition}                         \label{P1}
Each $Lh_0$-polar set with respect to $D$
is contained in a countable union of $N$-sets of the class 
$Lh_0(D)$.
\end{proposition}

Indeed, let $E$  be a  ${Lh_{0}}$-polar set with respect to $D$. Then by
Definition \ref{D2} there exists a function $\vartheta(x)\in Lh_0(D)$
such that $\vartheta(x){\not \equiv}-\infty$ and $\vartheta(x)|_E =
-\infty.$ Consider the sequence of functions
$\vartheta_k(x)=\frac1k\vartheta(x)$. Obviously, $\vartheta
_{k}(x)\in Lh_0(D)$ and $\vartheta _{k}(x){\not \equiv}-\infty$ and
$\vartheta _{k}(x)|_E=-\infty.$ Moreover, $\lim\limits_{k\rightarrow
\infty}\vartheta_k(x)=0$ for almost all $x\in D$ and
$\lim\limits_{k\rightarrow \infty}\vartheta_k(x)=-\infty$ for all
$x\in E$. It follows that
$$
E\subset \bigg\{x\in D:\lim \limits_{k\rightarrow \infty}
\vartheta_{k}(x)<\overline{\mathop {\lim} \limits_{y\rightarrow
x}}\lim \limits_{k\rightarrow\infty}\vartheta_k(x)\bigg\}.
$$
On the other hand, as we have showed above, the set
$$
\bigg\{x\in D:\lim \limits_{k\rightarrow \infty} \vartheta_{k}(x)<
\overline{\mathop {\lim} \limits_{y\rightarrow x}}\lim
\limits_{k\rightarrow\infty}\vartheta_k(x)\bigg\}
$$
consists of a countable union of $N$-sets of the class
$Lh_0(D)$.

\section
{\bf Proof of the theorem}

Suppose that the number $r$ is sufficiently large such that
$\widehat{V}\subset \mathbb{C}^m$ contains the closure of the unit
polydisc
$$
U=\bigg\{w=(w_1,w_2,...,w_m)\in
\mathbb{C}^m:|w_1|<1,|w_2|<1,...,|w_m|<1\bigg\},
$$
since the general case easily may be reduced to this case by linear
changing of $y$. By the lemma of \cite{sadullaev}, for each fixed $x\in
D$ the function $u(x,y)$ of variable $y$ can be extended holomorphly to
$\widehat{V}_r\subset \mathbb{C}^m$ and can be expanded into a  Hartogs'
series
$$
u(x,w)=\sum\limits_{|k|=0}^{\infty}c_k(x)w^k,
$$
where $k=(k_1,k_2,...k_m)$ is a multiindex. Clearly, the function
$u(x,y)$ is harmonic in each variable and coefficients $c_k\in
h(D)$. Also, according to the  Cauchy's inequalities (see
\cite{shabat}, \cite{vladimrov}), for each set $D_0:E\subset
D_0\subset \subset D$ the following estimation holds
$$
|c_k(x)|\leq M,~\forall k,~\forall x\in D_0,
$$
where $M=\sup\{|u(x,w)|:(x,w)\in \overline{D}_0\times
\overline{U}\}$, i.e., the sequence of functions $c_k(x)$ is locally
uniformly bounded in the domain $D.$ On the other hand, by the same lemma
from \cite{sadullaev}, for each fixed $x\in E$ the function $u(x,w)$
of the variable $w$ is holomorphic everywhere in $\mathbb{C}^m$, which means
that for all $x\in E$
 $$
 \overline{\mathop
{\lim}\limits_{|k|\rightarrow\infty}}\sqrt[|k|]{|c_k(x)|}=0.
 $$
 It follows that the sequence $\frac{1}{|k|}\ln |c_k(x)|$ is locally
uniformly bounded from above and
 $$
\overline{\mathop
{\lim}\limits_{|k|\rightarrow\infty}}\frac{1}{|k|}\ln
|c_k(x)|=-\infty,~\forall x\in E.
 $$
Put $\overline{\mathop {\lim}\limits_{y\rightarrow
x}}\,\,\,\,\overline{\mathop
{\lim}\limits_{|k|\rightarrow\infty}}\frac{1}{|k|}\ln
|c_k(y)|=\vartheta(x)$. Since the set
 $$
 F=\bigg\{x\in D:\overline{\mathop
{\lim}\limits_{|k|\rightarrow\infty}}\frac{1}{|k|}\ln
|c_k(x)|<\vartheta(x)\bigg\}
 $$
consists of the countable union of the $N$-sets of $ Lh_0(D)$
and $E$ is not contained in a countable union of the $N$-sets
of $ Lh_0(D)$, the set $E\setminus F=\{x\in
E:\vartheta(x)=-\infty\}$ is not contained in a countable
association of the $N$-sets of $Lh_0(D)$ as well. By Proposition \ref{P1} the
set $E\setminus F$ is not $Lh_0$-polar with respect to $D$.
Consequently, $\vartheta(x)\equiv-\infty$, i.e.,
$$
\overline{\mathop
{\lim}\limits_{|k|\rightarrow\infty}}\frac{1}{|k|}\ln
|c_k(x)|=-\infty,~\forall x\in D.
 $$
Thus, we get
 $$
\overline{\mathop
{\lim}\limits_{|k|\rightarrow\infty}}\sqrt[|k|]{|c_k(x)|}=0,~~\forall
x\in D.
 $$
Therefore, for all $x\in D$, the function  $u(x,y)$ of variable $y$ 
extends to a harmonic function on the whole $\mathbb{R}^m$. Hence,
$u(x,y)$ harmonically extends in $D\times \mathbb{R}^m$. The
theorem is proved.

{\em Authors' addresses}: {\em S.A. Imomkulov, Y.R. Saidov}, Urganch
State University, Urganch city, Khamid Alimjan -14, 740000,
Usbekistan.

e-mail: \texttt{ (S.A. Imomkulov), sevdi@rambler.ru, (Y.R. Saidov)
ysaidov@rambler.ru}

\end{document}